\documentclass[a4paper,11pt]{amsart}
\usepackage{amsmath,amsfonts,amssymb,epsf,artju}
\usepackage[dvipdfm]{hyperref}
\newcommand{\R}{\mathbb R}
\newcommand{\Z}{\mathbb Z}
\newcommand{\GD}{\Delta}
\newcommand{\GS}{\Sigma}

\newcommand{\p}{\partial}
\newcommand{\lk}{\operatorname{lk}}
\newcommand{\sminus}{\smallsetminus}

\theoremstyle{plain}

\newtheorem{thm}{\itshape\ Theorem}

\title{Linking number in
a projective space as the degree of a map}
\author{Julia Viro (Drobotukhina)}
\address{Department of Mathematics, Uppsala University, S-751 06
Uppsala, Sweden}
\email{julia@math.uu.se}
\begin{document}

\begin{abstract} For any two disjoint oriented circles embedded into the
3-dimensional real projective space, we construct a
3-dimensional configuration space and its map to the projective space
such that the linking number of the circles is the half of the degree of 
the map. Similar interpretations are given for the linking number of
cycles in a projective space of arbitrary odd dimension and the
self-linking number of a zero homologous knot in the 3-dimensional 
projective space.
\end{abstract}

\maketitle
\section*{Introduction}
As is well-known, the linking number of disjoint oriented circles $C_1,
C_2\subset\R^3$ can be presented as the degree of the map
\begin{equation}\label{*}
C_1\times C_2\to S^2: (x,y)\mapsto \frac{y-x}{|y-x|}.
\end{equation}
The linking number is defined in a more general situation for two  disjoint
oriented circles realizing homology classes of finite order in any
oriented 3-manifold. In particular, for disjoint oriented circles in
the projective space $\R P^3$. However, even in $\R P^3$, with its rich
geometry, the interpretation of the linking number as the degree of a map
similar to \eqref{*} does not exist.  The reasons for this are
presented below in Section \ref{s1}.

In this paper for any two disjoint oriented circles in $\R P^3$, we
construct  an oriented 3-dimensional configuration
space and its map to $\R P^3$ such that the degree of this map is the
linking number of the circles multiplied by 2. This construction seems
to be the closest possible replacement for \eqref{*}.  A diagramatical
formula for the linking number, which emerges in evaluation of this
degree by counting pre-images of a regular value, is similar to  the
well-known diagramatical formula for the degree of \eqref{*}. In fact,
a local ingredient in the diagramatical formula for the  linking number
of circles in $\R P^3$, the local writhe of a link diagram at a
crossing point, has appeared in my paper \cite{D1} and inspired the
present work.

The construction of the configuration space and its map to $\R P^3$
discussed above are generalized
to the case of a pair of disjoint oriented closed submanifolds (or even
sub-pseudo-manifolds) of $\R P^n$ such that the sum of their dimensions
equals $n-1$.

The first version of this paper appeared in 1998 as a preprint of
Mathematics Department Uppsala University (U.~U.~D.~M.\ Report 1998:29).
During my visit to MSRI in March-April 2004, I realized that the
result of that paper can be considered similar to recent results by
Jean-Yves Welschinger \cite{W1} and \cite{W2}. He considered real
versions of the following classical problem of enumerative geometry: 
how many rational curves of given degree pass through a generic 
collection of points. For an appropriate number of points, the number of
complex curves does not depend on the points, but the number of real
curves does depend. J.-Y.~Welschinger associated to each of the
real curves under consideration a sign in such a way that counting 
the curves with these signs gives a number which does not depend 
on the position of the points. 

The Main Theorem of this paper can be considered as a similar result
concerning a similar problem of enumerative geometry: counting of 
real lines in $\R P^3$ which intersect each of
two given oriented closed curves and pass through a given point. 
Indeed, if each of the lines is taken with the sign equal to the 
local degree of the map $F$ (see Section \ref{s4} below) at the 
corresponding point  
(this local degree is the writhe number of the corresponding crossing
point in the projection of the curves from  the point), 
then the result of count 
does not depend on the point and equals the linking number of the curves
multiplied by two.

I am grateful to Mathematics Department of Uppsala University and
MSRI for a partial support of this work.

\section{The classical construction and impossibility of its
straightforward generalization}\label{s1}

The Gauss map \eqref{*} can be viewed as a composition of the natural
inclusion
$$C_1\times C_2\to 
\{(x,y)\in\R^3\times\R^3\mid x\ne y\}=(\R^3\times\R^3)\sminus\Delta $$
and the projection
$$(\R^3\times\R^3)\sminus\Delta\to S^2 :
 (x,y)\mapsto \frac{y-x}{|y-x|}.$$
This composition
induces a homomorphism
$$\Z=H_2(C_1\times C_2)\to H_2(\R^3\times\R^3\sminus\Delta))\to
H_2(S^2)=\Z.$$
The degree of the Gauss map is an integer $deg$ such that
$[C_1\times C_2]\mapsto deg \cdot[S^2]$ for fundamental classes  $[C_1\times
C_2]$ and $[S^2]$. Since
$$H_2(\R^3\times\R^3\sminus\Delta)\to H_2(S^2)=\Z$$
is an isomorphism as a map induced by  a homotopy
equivalence, the degree of the Gauss map is actually
 defined by the inclusion homomorphism $ H_2(C_1\times C_2)
\to H_2(\R^3\times\R^3\sminus\Delta)$.

Consider now the case of two oriented disjoint circles $C_1, C_2$ in
$\R P^3$. The classical case suggests to consider the inclusion
$$C_1\times C_2\to \{(x,y)\in\R P^3\times\R P^3 : x\ne y\}=
(\R P^3\times\R P^3)\sminus\Delta .$$
A natural hope is to find an expression for the linking number in terms
of the homomorphism induced by this map in the second homology.
However, this hope is ruined by the fact that the second homology of
the target space is $\Z_2$.

To prove that $H_2((\R P^3\times\R P^3)\sminus\Delta)=\Z_2$, consider
the natural projection  $$(\R P^3\times\R P^3)\sminus\Delta\to\R P^3,
(x,y)\mapsto x.$$ This is a locally trivial fibration over $\R P^3$ with
fiber  $\R P^3\sminus\{x\}$. The fiber is homotopy equivalent to 
$\R P^2$ and the fibration can be trivialized using a trivialization of
the tangent bundle of $\R P^3$.  
Hence $H_2((\R P^3\times\R
P^3)\sminus\Delta)=H_2(\R P^3\times\R P^2)$.  By K\"unneth formula,
$ H_2(\R P^3\times\R P^2)=\Z_2$.

Thus the inclusion $C_1\times C_2\to (\R P^3\times\R P^3)
\sminus\Delta$ induces in the second homology a homomorphism $\Z\to\Z_2$
and the linking number of $C_1$ and $C_2$ (which is a half-integer) cannot
be expressed in terms of this homomorphism.

\section{Configuration space}\label{s2}
Let  $C_1, C_2$ be 
two oriented disjoint circles in $\R P^3$. Consider a space $K_0$ of
all triples $(x,t,y)$, where $x\in C_1$, $y\in C_2$ and $t$
($x\ne t\ne y $) is a point on the line passing through $x$ and $y$.

This space is naturally equipped with a structure of a smooth oriented
3-manifold. We introduce this structure by describing an appropriate
class of local coordinate systems. One of the coordinates is defined
globally in the following way.
For $(x,t,y)\in K_0$, consider a straight  segment $I$  containing $t$ and
bounded by $x,y$. Denote by $J$ the segment on $I$ bounded by $x$ and
$t$. Put $\tau=\frac{|J|}{|I|}$, where $|\cdot|$ denotes the
length with respect to the Fubini-Study metric on $\R P^3$. This
defines a function $\tau:K_0\to(0,1)$.

To construct a local coordinate system in
a neighborhood of an arbitrary point $(x_0,t_0,y_0)$ of $K_0$, choose
 a  local coordinate $\xi$ in a neighborhood of $x_0$ on $C_1$,  a
local coordinate $\eta$ in a neighborhood of $y_0$ on $C_2$.
Functions $\xi$, $\tau$, $\eta$ restricted to a neighborhood of
$(x_0,t_0,y_0)$ form a local coordinate system in $K_0$.

The coordinate systems provided by this construction define a differential
structure on $K_0$. By choosing only those of these coordinate systems,
in which the first coordinate defines the orientation of $C_1$
and the third coordinate, the orientation of $C_2$, we define an orientation on
$K_0$. This orientation is canonical in the sense that it depends only on
the orientation of $C_1$ and $C_2$ and the ordering of them (the
orientation changes, if we exchange $C_1$ and $C_2$). 

In order to make this configuration space  compact, we allow $t=x$ or
$t=y$. In other words, consider
$$K_1=\{(x,t,y) \in C_1\times \R P^3\times C_2 \mid t\in Span(x,y)\}.$$
This is a compact smooth 3-manifold containing $K_0$ as the complement
of two disjoint copies of $C_1\times C_2$. These copies are
$S_1=\{(x,t,y)\in K_1 \mid t=x\}$ and $S_2=\{(x,t,y)\in K_1\mid t=y\}$.
Denote $S_1\cup S_2$ by $S$.

The orientation of $K_0$ does not extend over $K_1$: it reverses when
one goes transversally through the surface $S$ at any point. Moreover,
if at least one of the circles $C_1$ and $C_2$ is not homologous to
zero in $\R P^3$, then $K_1$ is not orientable. Indeed, $K_1$ is
fibered over the torus $C_1\times C_2$ with fiber $S^1$,  and the first
Stiefel-Whitney class is zero iff both $C_1$ and $C_2$ are
zero-homologous in $\R P^3$. (To describe the fibration completely,
observe that the Euler number of this fibration is 0, since it has
disjoint sections $S_1$ and $S_2$.)

\section{Mapping and final factorizing of the configuration
space}\label{s3}

There is a natural map $F_1:K_1\to \R P^3$ defined by $F(x,t,y)=t$.
Under this map the whole surface $S_i$ is mapped to curve $C_i$,
$i=1,2$. The preimage of a point $x\in C_1$ is the circle formed by
triples $(x,x,y)$ with $y\in C_2$, while the preimage of $y\in C_2$
consists of $(x,y,y)$ with $x\in C_1$.

Let us contract each of these circles. In other words, consider the
quotient space
$$K=K_1/_{
(x_1,y,y)\sim(x_2,y,y),\quad(x,x,y_1)\sim(x,x,y_2)}.$$
This space consists of an intact copy of $K_0$, and the images of $S_1$ and
$S_2$, which are naturally identified with $C_1$ and $C_2$, respectively.

One can show that \begin{itemize}
\item If both $C_1$ and $C_2$ are not zero-homologous in
$\R P^3$ then $K$ is homeomorphic to $\R P^3$.
\item If one of $C_1$ and
$C_2$ is zero-homologous, while the other is not, then $K$ is
homeomorphic to the quotient space of $S^3$ which can be obtained by
factorizing of an unknot in $S^3$ by a fixed point free involution.
\item If both $C_1$ and $C_2$ are zero-homologous in $\R P^3$, then $K$
is homeomorphic to the result of attaching to each other two copies of
the 3-sphere $C_1*C_2$ by the identity homeomorphism between the copies
of $C_1\cup C_2\subset C_1*C_2$.\end{itemize}
We do not prove this, since it is not needed in what follows.

Thus $K$ is not necessarily a 3-manifold, but it is a stratified space
which consists of a 3-stratum $K_0$ and 1-strata $C_1$ and $C_2$.
Hence the orientation of $K_0$ defines a homology class $[K_0]$, which
belongs to $H_3(K,C_1\cup C_2)$. The latter group is naturally isomorphic
to $H_3(K)$: the inclusion map $H_3(K)\to H_3(K,C_1\cup C_2)$ is an
isomorphism, since it is surrounded in the homology sequence of the
pair $(K,C_1\cup C_2)$ by trivial groups $H_3(C_1\cup C_2)$ and
$H_2(C_1\cup C_2)$. Therefore we identify $H_3(K,C_1\cup C_2)$ with
$H_3(K)$. Denote the homology class corresponding to $[K_0]$ under this
identification by $[K]$.

\section{The Main Theorem}\label{s4}

The map $F_1:K_1\to \R P^3$ defined above defines a map $F:K\to\R P^3$.
The induced map $F_*:H_3(K)\to H_3(\R P^3)(=\Z)$ maps $[K]$ to a class
which can be presented as $d[\R P^3]$, since the orientation class
$[\R P^3]$ of $\R P^3$ generates $H_3(\R P^3)$. The integer coefficient 
$d$ depends only on the initial data of the construction which gave rise
to $K$ and $F$. That is it depends only on $C_1$ and $C_2$. 
Let us denote $d$ by $d(C_1,C_2)$. 

\begin{thm}\label{th1}
Under the conditions above, $d(C_1,C_2)=2\lk(C_1,C_2)$.
\end{thm}

\begin{proof} To evaluate $d(C_1,C_2)$ geometrically, we choose a
regular value $v\in \R P^3\sminus(C_1\cup C_2)$ of the map $F$ and take
the sum of local degrees of $F$ over the preimages of $v$.

We have to prove that this sum is equal to the doubled linking number
of $C_1$ and $C_2$. To evaluate the doubled linking number, we take the
intersection number of $C_2$ with a chain $\GS$ such that $\p\GS=2C_1$.
For $\GS$, we choose the projective cone $C_v(C_1)$ over $C_1$ with
vertex at $v$, i.e. the union of all projective lines passing through
$v$ and intersecting $C_1$. Obviously, points of $C_2\cap C_v(C_1)$ are
in a natural one-to-one correspondents with points of $F^{-1}(v)$, and
we need to check that the local intersection number of $C_2$ and
$C_v(C_1)$ at a point of $C_2\cap C_v(C_1)$  is equal to the local
degree of $F$ at the corresponding point of $F^{-1}(v)$. The comparison
of these two local numbers is shown in Figure \ref{loc-or}.
\begin{figure}[hbt]
\centerline{\epsffile{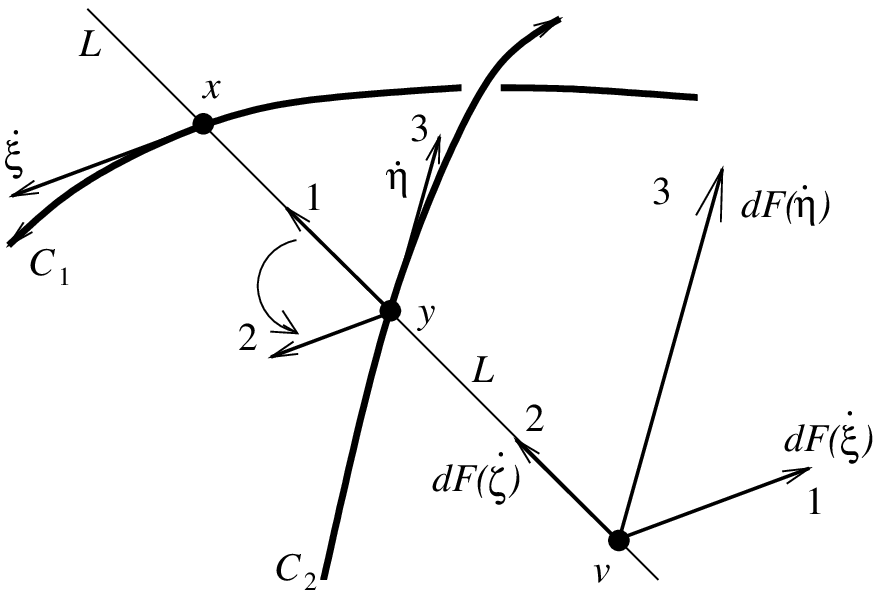}}
\caption{}
\label{loc-or}
\end{figure}
There we see a line $L$ intersecting oriented embedded circles $C_1$ and
$C_2$ in points $x$ and $y$, respectively. Vectors 
$\dot\xi=\frac{\p}{\p\xi}$ and $\dot\eta=\frac{\p}{\p\eta}$ are tangent
to $C_1$ and $C_2$ at $x$ and $y$, respectively.

By the definition of differential, $dF(\dot\xi)$ is the velocity of $t$
at $t=v$ as $x$ moves along $C_1$ with velocity $\dot\xi$, while $y$ 
and the ratio, in which $t$ divides $[y,x]$, stay fixed (the ratio is
equal to the ratio, in which $v$ divides $[y,x]$).  Therefore $L$,
$\dot\xi$, $dF(\dot\xi)$ are coplanar and $\dot\xi$ is contained in $L$
iff  $dF(\dot\xi)$ is contained in $L$. Similarly, $dF(\dot\eta)$ is
the velocity of $t$ at $t=v$ as $y$ moves along $C_2$ with velocity
$\dot\eta$, while $x$  and the ratio in which $t$ divides $[y,x]$ stay
fixed.   Therefore $L$, $\dot\eta$, $dF(\dot\eta)$ are coplanar and
$\dot\eta$ is contained in $L$ iff  $dF(\dot\eta)$ is contained in $L$.

Therefore, if $L$, $\dot\xi$ and $\dot\eta$ are not coplanar, then
$dF$ is not degenerate at $t$. The local degree of $F$ at $(x,t,y)$ is
the value of the standard orientation of $\R P^3$ on the frame
$dF(\dot\xi), dF(\dot\tau),dF(\dot\eta)$. (Recall that we defined the
orientation of $K_0$ as positive on the local coordinate system
$\xi,\tau,\eta$, see Section \ref{s2}.)

Let us compare this local degree to the local intersection number of
$C_2$ and $C_t(C_1)$ at $y$. The local intersection number is the value
taken by the same standard orientation of $\R P^3$ on the frame which
consists of $\dot\eta$ and the basis of the tangent space of $C_t(C_1)$
at $y$ positively oriented with respect to the orientation such that
$\p C_t(C_1)=2C_1$. One of frames of this sort is shown in Figure
\ref{loc-or}. Its first vector is directed along $L$  outwards the
segment $[t,y]$ which does not contain $x$. The second vector is
directed to the same side of the segment $[x,y]$ not containing $t$ as
$\dot\xi$.

In Figure \ref{loc-or} the orientations of these frames obviously coincide.
(For reader's convenience the orderings of bases vectors are shown by
numbers.)
Figure \ref{loc-or} is necessarily affine, and, at first glance, quite
special. One could guess that it does not represent all possible
configurations. However, in fact, it is general: we get affine picture
shown in Figure \ref{loc-or} for any $(x,t,y)\in K_0$, if we choose the
infinity plane intersecting  $L=span(x,y)$ in the segment bounded
by $x$ and $y$ and not containing $t$.  \end{proof}

\section{Versions and generalizations of Main Theorem}\label{s5}

The number $d(C_1,C_2)$ involved in Theorem \ref{th1} can be introduced in a
different way.   Instead of $F$ and
$K$ one can use $F_1$ and $K_1$. This would eliminate a step of the
construction, which was the main contents of Section \ref{s3}.
However then instead of $[K]\in H_3(K)$ one should use relative
homology $H_3(K_1,S)$ and the homology class $[K_1]\in H_3(K_1,S)$
defined by the orientation of $K_0(=K_1\sminus S)$. The role of
$F_*:H_3(K)\to H_3(\R P^3)$ would be played by $F_{1*}:H_3(K_1,S)\to
H_3(\R P^3,C_1\cup C_2)$. Since  $H_3(\R P^3,C_1\cup C_2)$ is
isomorphic to $H_3(\R P^3)=\Z$, the image of $[K_1]$ under $F_{1*}$ is
equal to the $d$-fold multiple of the generator of
$H_3(\R P^3,C_1\cup C_2)$ defined by the standard orientation of $\R
P^3$. Of course, this $d$ is equal to $d(C_1,C_2)$.

Theorem \ref{th1} can be generalized by admitting more general $C_1$
and $C_2$. In the original setup, $C_1$ and $C_2$ are oriented smooth
disjoint circles embedded in $\R P^3$. Of course, the assumption of
their connectedness is not needed. If $C_1$, $C_2$ are oriented
disjoint closed 1-submanifolds of $\R P^3$, one can repeat everything
made above, besides the discussion of the topological type of $K$ in
Section \ref{s3}. This discussion can be skipped, it was included only
to show that $K$ may have singularities. In the case of multicomponent
$C_1$ and $C_2$, $K_1$ consists of connected components which are of
the same types as above, while $K$ can be obtained from pieces of the
same type by gluing them together along the components of $C_1$ and
$C_2$. The formulation and proof of Theorem \ref{th1} do not change.

Moreover, Theorem \ref{th1} can be generalized to the case of cycles
$C_1$, $C_2$ supported by disjoint smoothly stratified 1-dimensional
spaces $|C_1|$ and $|C_2|$ (i.e., graphs). Each 1-stratum of $|C_i|$ is
oriented and equipped with a coefficient so that the corresponding
linear combination of the orientation cycles of the strata is $C_i$.
Then we define
$$K_1=\{(x,t,y)\in|C_1|\times\R P^3\times|C_2|\mid t\in Span(x,y)\},$$
$S_1=\{(x,t,y)\in K_1\mid t=x\}$ and $S_2=\{(x,t,y)\in K_1\mid t=y\}$,
put $S=S_1\cup S_2$ and $K_0=K_1\sminus S$. The stratifications of
$|C_1|$ and $|C_2|$ define natural stratifications of these spaces.
The orientations of 1-strata of $|C_i|$ define orientations of 3-strata
of $K_0$ as above. Each of these 3-strata is defined by a pair
consisting of a 1-stratum of $|C_1|$ and a 1-stratum of $|C_2|$.
Assigning to each of the 3-strata the products of the corresponding
coefficients, we get a 3-cycle of
$$K=K_1/_{(x_1,y,y)\sim(x_2,y,y) \quad (x,x,y_1)\sim(x,x,y_2)}.$$
The rest of the preparations to Theorem \ref{th1}, the theorem itself
and its proof can be repeated literally.

Consider now a straightforward high-dimensional generalization.
It is possible  to consider arbitrary cycles, as above, but we restrict
ourselves to the case of submanifolds.
Let $C_1$ and $C_2$ be oriented closed disjoint smooth submanifolds of
dimensions $p_1$ and $p_2$, respectively, of the real projective space
$\R P^{n}$. Assume that $n=p_1+p_2+1$ is an odd number and, if $p_i=0$,
then $C_i$ is zero-homologous in $\R P^n$ (i.e., $C_i$ consists of an
even number of points, half of which are oriented positively and half,
negatively).

As is well known, in this situation there is an integer or half-integer
$\lk(C_1,C_2)$, which can be defined as follows. Take any oriented
$(p_1+1)$-dimensional chain $\GS$ with integer coefficients in $\R P^n$
with $\p\GS=2C_1$. If $p_1>0$, one can construct such $\GS$ as a
projective cone, see Proof of Theorem \ref{th1}.
Put it to general position with respect to $C_2$.
Then $\lk(C_1,C_2)$ is $\frac12\GS\circ C_2$, where $\GS\circ C_2$ is
the intersection number.

Define
$$K_1=\{(x,t,y)\in C_1\times\R P^n \times C_2\mid t\in Span(x,y)\},$$
$S_1=\{(x,t,y)\in K_1\mid t=x\}$ and $S_2=\{(x,t,y)\in K_1\mid t=y\}$,
put $S=S_1\cup S_2$ and $K_0=K_1\sminus S$. The smooth structures and
orientations of $C_1$ and $C_2$ define a smooth structure and
orientation of $K_0$ exactly as in Section \ref{s2}. In particular, the
orientation is positive on the local coordinate system, in which the
first $p_1$ coordinates are taken from the positively oriented local
coordinates on $C_1$, the $p_1+1^{\text{st}}$ coordinate is directed
along the line connecting $C_1$ to $C_2$ towards $C_2$, and the last
$p_2$ coordinates are taken from the positively oriented local
coordinates on $C_2$. The natural map
$$F_1:K_1\to\R P^n:(x,t,y)\mapsto t$$
can be factored through the quotient space
$$K=K_1/_{(x_1,y,y)\sim(x_2,y,y) \quad (x,x,y_1)\sim(x,x,y_2)}.$$
Denote, as above, the resulting map $K\to \R P^n$ by $F$. The
orientation of $K_0$ defines a class $[K]\in H_n(K)$. The image
$F_*[K]\in H_n(\R P^n)$ is $d(C_1,C_2)[\R P^n]$.
As in Theorem \ref{th1}, under
these conditions $d(C_1,C_2)=2\lk(C_1,C_2)$.

\section{The self-linking coefficient}\label{s6}

In \cite{D1}, I introduced a numerical invariant of a homologous to zero
oriented knot in $\R P^3$. The preimage of such a knot in the universal
covering space $S^3$ of $\R P^3$ is a two-component link. The
{\it self-linking coefficient\/} is defined to be the linking number of the
components of this link. The self-linking coefficient was used in
\cite{D1} as an obvious obstruction for a knot to be isotopic to a knot
contained in an affine part of $\R P^3$.

Below the self-linking coefficient $sl(k)$ of a knot $k$ is represented
via the degree of a map. The source space of the map is constructed as
follows.

For a triple $(x,t,y)$ consisting of $x,y\in k$, $x\ne y$ and 
$t\in Span(x,y)$,
consider a loop composed of a half of $k$, which is bounded by $x$ and
$y$, and the segment of the line $Span(x,y)$ which is bounded by $x$ and $y$
and does not contain $t$. The homology class of this loop does not
depend on the choice of the half of $k$, since $k$ is zero-homologous
in $\R P^3$. Denote by $K_0$ the set of triples $(x,t,y)$ with
$x,y\in k$, $x\ne y$ and $t\in Span(x,y)$ such that the homology class of the
loop constructed above for $(x,t,y)$ is not zero.

The orientation of $k$ and the direction of the segment of $Span(x,y)$
from $x$ to $y$ define an orientation of $K_0$.

The space $K_0$ is fibered over  $k\times k\sminus\GD$ with fiber $\R$.
Consider a larger space $K_1$ which is the closure of $K_0$ in the
space of all triples $(x,t,y)$ with $x,y\in k$ and $t\in Span(x,y)$.
Let $K$ be the quotient space of $K_1$:
$$K=K_1/_{(x,x,y_1)\sim(x,x,y_2)\quad (x_1,y,y)\sim(x_2,y,y)}.$$
Obviously, $K\sminus K_0$ can be identified with $k$. Hence $K$ is a
stratified space with a 3-dimensional stratum $K_0$ and 1-dimensional
stratum $k$. The orientation of $K_0$ defines a homology class
$[K]\in H_3(K)$.

Define $F_1:K_1\to\R P^3$ by $F_1(x,t,y)=t$. Denote by $F$ the quotient
map $K\to \R P^3$. Let $d(k)$ be an integer such that $F_*[K]=d(k)[\R
P^3]$.

\begin{thm}\label{th2}
Under the conditions above
$$d(k)=2sl(k).$$
\end{thm}

The proof of Theorem \ref{th2} is similar to the proof of Theorem
\ref{th1}. \qed

\end{document}